\documentclass[11pt]{amsart}

\usepackage{url}
\usepackage{graphicx}
\usepackage{amssymb, amscd, amsmath}
\usepackage{xy}
\usepackage{bm}
\xyoption{all}
\usepackage[linktocpage]{hyperref}
\usepackage[dvipsnames]{xcolor}
\usepackage{outline}
\usepackage{epstopdf}
\usepackage{rotating}
\usepackage{mathabx}
\usepackage{fancyhdr}
\usepackage[left=1.25in,top=1.2in,right=1.25in,bottom=1.2in,head=.2in]{geometry}
\usepackage{mathabx}

\usepackage{pinlabel}

\newtheoremstyle{newstyle}% name
  {.8\baselineskip\@plus.2\baselineskip\@minus.2\baselineskip}% Space above
  {.4\baselineskip\@plus.2\baselineskip\@minus.2\baselineskip}% Space below
  {\slshape}% Body font
  {}%Indent amount (empty = no indent, \parindent = para indent)
  {\bfseries}%  Thm head font
  {.}%       Punctuation after thm head
  { }%      Space after thm head: " " = normal interword space;
  {}%       Thm head spec

\theoremstyle{newstyle}
\newtheorem{thm}{Theorem} % lettered theorems (A,B,C,D)

\newtheorem{cor}[thm]{Corollary}

\newtheorem{theorem}{Theorem}[section]
\newtheorem*{theorem*}{Theorem}
\newtheorem{lemma}[theorem]{Lemma}

\newtheorem*{corollary*}{Corollary}
\theoremstyle{definition}
\newtheorem{definition}[theorem]{Definition}

\newtheorem*{remark*}{Remark}
\newtheorem*{remarks*}{Remarks}
\newtheorem*{addenda*}{Addenda}

\newcommand{\fig}[3]{\begin{figure}[h!] \includegraphics[height=#1pt]{#2}#3\end{figure}}

\newcommand{\figref}[1]{Figure~\ref{F:#1}}
\newcommand{\secref}[1]{Section~\ref{S:#1}}
\newcommand{\thmref}[1]{Theorem~\ref{T:#1}}
\newcommand{\lemref}[1]{Lemma~\ref{L:#1}}

\newcommand{\bz}{\mathbb Z}

\newcommand{\bc}{\mathbb C}

\newcommand{\bu}{\mathbb U}
\newcommand{\bv}{\mathbb V}
\newcommand{\bw}{\mathbb W}

\newcommand{\calm}{\mathcal M}
\newcommand{\caln}{\mathcal N}

\newcommand{\cs}{\mathbin{\#}} % or {\,{\hash}\,}

\newcommand{\st}{\,\vert\,}

\newcommand{\hookto}{\hookrightarrow}

\newcommand{\id}{\textup{id}}

\newcommand{\cponebar}{\overline{\bc\textup{P}}\,\!^1}
\newcommand{\cptwo}{\bc\textup{P}^2}
\newcommand{\cptwobar}{\overline{\bc\textup{P}}\,\!^2}
\DeclareMathOperator{\SO}{SO}
\DeclareMathOperator{\interior}{int}

\newcommand{\G}{G}

\def\plus{\natural\,}
\def\sum{\sharp\,}

\newcommand{\Z}{\mathbb Z}

\newcommand{\R}{\mathbb R}

\title{Equivariant corks}
\author[Dave Auckly]{Dave Auckly${}^1$}
  \address{Department of Mathematics
  \newline\indent Kansas State University
  \newline\indent  Manhattan, Kansas 66506}
  \email{dav@math.ksu.edu}
\author[Hee Jung Kim]{Hee Jung Kim${}^{1,2}$}
  \address{Department of Mathematical Sciences 
  \newline\indent Seoul National University
  \newline\indent Seoul 151-747, South Korea}
  \email{heejungorama@gmail.com}
\author[Paul Melvin]{Paul Melvin${}^1$}
  \address{Department of Mathematics
  \newline\indent Bryn Mawr College
  \newline\indent  Bryn Mawr, PA 19010}
  \email{pmelvin@brynmawr.edu}
\author[Daniel Ruberman]{Daniel Ruberman${}^{1,3}$}
  \address{Department of Mathematics, MS 050
  \newline\indent Brandeis University 
  \newline\indent Waltham, MA 02454}
  \email{ruberman@brandeis.edu}

\thanks{Partially supported by ${}^1$AIM SQuaRE grant, ${}^2$NRF grants 2015R1D1A1A01059318 and BK21 PLUS SNU Mathematical Sciences Division, ${}^3$NSF Grants 1105234 and 1506328, and NSF FRG Grant 1065827.
}

\begin{document}
\begin{abstract} For suitable finite groups $\G$, we construct contractible $4$-manifolds $C$ with an effective $\G$-action on $\partial C$ whose associated pairs $(C,g)$ for all $g \in \G$ are distinct smoothings of the pair $(C,\partial C)$.  Indeed $C$ embeds in a $4$-manifold so that cutting out $C$ and regluing using distinct elements of $\G$ yield distinct smooth $4$-manifolds.  
\end{abstract}
\maketitle
\vskip-.3in
\vskip-.3in

\setcounter{section}{-1}

%%%%%%%%%%%%%%%%%%%%%%%%%%%%%%%%%%%%%%%%%%%%%%%%%%%%%%%%%%%%%%%%%%
%%%%%%%%%%%%%%%%   SECTION 0 : INTRODUCTION  %%%%%%%%%%%%%%%%%%%%%
%%%%%%%%%%%%%%%%%%%%%%%%%%%%%%%%%%%%%%%%%%%%%%%%%%%%%%%%%%%%%%%%%%

\section{Introduction}\label{S:intro}

\def\u#1{\underline{#1}}

A {\it cork} is a smooth, compact, contractible $4$-manifold with an involution on its boundary that does not extend to a diffeomorphism of the full manifold.   Akbulut\,\cite{akbulut:contractible} discovered this phenomenon for the classical Mazur manifold $\bw$\,\cite{mazur:contractible} with the boundary involution $\tau$ shown in Figure 1,
%\figref{MazurCork}, 
proving that $\bw$ embeds in a $4$-manifold $X$ so that the result of removing $\bw$ and regluing it using $\tau$ is not diffeomorphic to $X$.  

%%%%%%%%%% FIG 1 %%%%%%%%%%
\fig{80}{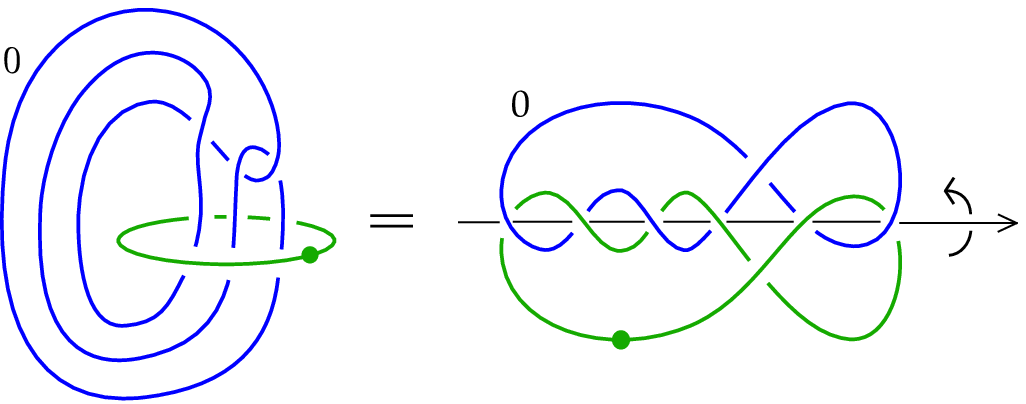}{
  \put(-10,45){\footnotesize$\tau$}
\caption{The Mazur Cork}
\label{F:MazurCork}}
%%%%%%%%%%%%%%%%%%%%%%%%%%%

\noindent This operation is called {\it cork twisting}, and it is now known ~\cite{curtis-freedman-hsiang-stong,matveyev:h-cobordism} that any two smooth, closed, simply-connected $4$-manifolds that are homeomorphic differ by a single cork twist.  It is not known whether the same cork can be used in all situations, i.e.\ whether there exists a {\it universal cork}\/; indeed it is conceivable, though unlikely, that the Mazur cork is universal.
 
The property that the cork twist $\tau$ is an {\it involution} is interesting, indeed inherent in most constructions of corks to date, but it is not clear that it is fundamental to the relation between cork twists and other smooth $4$-manifold constructions. It is therefore natural to ask whether cutting and gluing by higher order diffeomorphisms of the boundary of a contractible submanifold of a $4$-manifold can change the underlying smooth structure.  In this note we give an affirmative answer, producing examples of embeddings of contractible $4$-manifolds with twists of arbitrary finite order that alter the ambient smooth structure; a different construction of such examples was given in a recent preprint of Tange \cite{tange:cycliccorks}. 

In fact we show more:  For suitable finite groups $G$, there exist contractible $4$-manifolds with effective $G$-actions on the boundary that embed in closed $4$-manifolds so that twists corresponding to distinct elements of $G$ yield distinct smooth structures.  We call such a gadget an {\it equivariant cork}, or $G$-cork if we want to specify the group.  

%%%%%%%% Theorem A %%%%%%%%

\begin{thm}\label{T:equiv-cork}
There exist $G$-corks for any finite subgroup $G$ of $\SO(4)$.
\end{thm}

If the action of $G$ on $S^3$ is free, then the action of $G$ on the boundary of cork constructed in the theorem is free; this seems to be a new phenomenon, even for $G= \Z_2$.  In a final section, we extend the notion of $G$-cork to the setting where the group $G$ is a subgroup of the mapping class group of the boundary; see the end of \secref{equiv} for details. We call this extension a {\em weakly equivariant} cork, and give an example of an effective weak $G$-cork in this sense where $G$ is a group that does not act effectively on any homology $3$-sphere.

\begin{thm}\label{T:w-eq}
There are groups $G$ that do not act effectively on any homology sphere, but for which there exist weakly equivariant $G$-corks.
\end{thm}

The boundaries of the corks constructed in the proof of \thmref{equiv-cork} are reducible.
%$3\text{-manifolds}$.  
In a sequel we will prove the following theorem, using rather different techniques from those in the current paper.

%%%%%%% Theorem B %%%%%%%

\begin{thm}\label{hypthm}
Given an oriented $3$-manifold $Y$ with an effective, orientation-preserving, smooth action of a finite group $\G$, there is an equivariant invertible $\mathbb{Z}[\pi_1(Y)]$-homology cobordism from it to a hyperbolic manifold. 
\end{thm}

As in~\cite{akbulut-ruberman:absolute}, this immediately implies:

%%%%%%% Corollary C %%%%%%

\begin{cor} \label{T:hyp-cork}
\ For any given finite subgroup \,$G$\, of $\SO(4)$, there exists a $G$-cork with hyperbolic boundary.
\end{cor}

Some experimentation with Snappy~\cite{SnapPy} suggests that the simplest $\Z_n$-corks in Tange's paper have hyperbolic boundaries, but a proof in general would require different techniques.

\subsection*{Acknowledgements}
The construction of G-corks and their hyperbolization were worked out by the authors at the American Institute of Mathematics (AIM) at our SQuaRE meeting in July 2015.  We thank AIM for its support for this and future endeavors. Our results were announced at the 2016 Joint Mathematics Meeting; see~\cite{akmr-abstract}.

%%%%%%%%%%%%%%%%%%%%%%%%%%%%%%%%%%%%%%%%%%%%%%%%%%%%%%%%%%%%%%%%%%
%%%%  SECTION 1 : PRELIMINARIES AND STATEMENT OF RESULTS  %%%%%%%%
%%%%%%%%%%%%%%%%%%%%%%%%%%%%%%%%%%%%%%%%%%%%%%%%%%%%%%%%%%%%%%%%%%

\section{Preliminaries and Statement of Results}\label{S:equiv}

In this section, we lay the groundwork for our proof of the existence of equivariant corks.  Most of the ideas discussed here are well known, but since we will use ``corks" in a broader sense than usual, and employ cork twists on multiple copies of boundary sums of embedded copies of the Mazur cork, we must give careful definitions of the relevant notions. 

\vskip-.05in
\vskip-.05in
\subsection*{Corks and Boundary Equivalence}

Extending the usual terminology, a {\it cork} will refer to any pair $(C,g)$ where $C$ is a smooth, compact, contractible $4$-manifold, and $g$ is an {\sl arbitrary} diffeomorphism of $\partial C$.  In particular $g$ need not be an involution, nor even of finite order, and $C$ need not be Stein (as is often assumed, cf.\,\cite{akbulut-yasui:corks-plugs}).  But if $g$ {\it is} a {\it special involution} (meaning orientation preserving with nonempty fixed point set, as with the Mazur twist $\tau$) then we also refer to $(C,g)$ as a {\it special $2$-cork}.   

In general, we call a cork $(C,g)$ {\it trivial} if $g$ extends to a diffeomorphism of $C$ (it always extends to a homeomorphism by \cite{freedman:simply-connected}) and {\it nontrivial} otherwise; with this convention, $(B^4,g)$ is a trivial cork for any $g$, whereas the Mazur cork $(\bw,\tau)$ is nontrivial.  These notions induce an equivalence relation on corks associated with the same underlying manifold: $(C,g)$ and $(C,h)$ are {\it boundary equivalent} if and only if $(C,g^{-1}h)$ is trivial, i.e.\ $g^{-1}h$ extends over $C$.

\subsection*{Boundary Sums of Corks} 

The {\it boundary sum} operation $\plus$ is well-defined on boundary equivalence classes of corks, as follows:  Given corks $(C_1,g_1)$ and $(C_2,g_2)$, choose (for $i=1,2$) diffeomorphisms $h_i$ isotopic (and thus boundary equivalent) to $g_i$ that are the identity on 3-balls $B_i\subset \partial C_i$.  Form $C_1 \plus C_2$ by identifying the $C_i$'s along the $B_i$'s so that $h_1$ and $h_2$ glue together to form $h_1\sum h_2$.  The result 
$$
(C_1,g_1)\,\plus(C_2,g_2) \ := \ (C_1 \plus C_2,h_1 \sum h_2)
$$
may depend on the choices of $h_i$ and $B_i$, but its boundary equivalence class does not.  Note however that $\plus$ {\it is} well-defined for special $2$-corks {\it without} imposing boundary equivalence; just choose the $B_i$ to be $g_i$-invariant $3$-balls centered at fixed points, and then $g_1\sum g_2$ is a well-defined involution, independent of the choices up to equivariant diffeomorphism.  

\vskip -5pt
\vskip -5pt
\subsection*{Cork Embeddings}   

A {\it cork embedding} of $(C,g)$ in a 4-manifold $X$ is a smooth embedding $e\!:\!C\hookto X$ with the induced map $\bar g=ege^{-1}$ on the boundary of its image $\overline C = e(C)$.  The associated {\em cork twist} $X_g^e$ is obtained by removing $\overline C$ from $X$ and regluing using $\bar g$\,:
$$
X_g^e \ = \ (X - \interior\overline C) \,\cup_{\bar g}\, \overline C\,.
$$
The embedding $e$ is {\it trivial} if $X_g^e$ is diffeomorphic to $X$, and {\it nontrivial} or {\it effective} otherwise.  Thus the nontriviality of $(C,g)$ can be verified by producing a nontrivial embedding, rather than trying to show directly that $g$ does not extend smoothly across $C$.

Note that the definition of boundary equivalence of cork maps is compatible with the use of such maps in changing smooth structures, because the result of twisting by $g$ is the same as the result of twisting by $h$ when $g^{-1}h$ extends across $C$. Conversely, given any nontrivial cork $(C,g)$, Akbulut and Ruberman construct a pair of absolutely exotic structures on a contractible manifold related by twisting $(C,g)$ \cite{akbulut-ruberman:absolute}. It follows that for any two boundary inequivalent diffeomorphisms $g$ and $h$, there is a $4$-manifold $X$ and an embedding $e:C\hookto X$ so that $X_g^e$ is not diffeomorphic to $X_h^e$.   Akbulut has made a similar observation.

\vskip -5pt
\vskip -5pt
\subsection*{Boundary Sums of Cork Embeddings}   

Given any pair of embeddings $e_i:C_i\hookto X$ (for $i=1,2$) of corks $(C_i,g_i)$ with disjoint images $\overline C_i=e_i(C_i)$ and induced boundary maps $\bar g_i:\partial \overline C_i\to\partial \overline C_i$, both twists can be performed simultaneously to produce the 4-manifold 
$$
X_{g_1g_2}^{e_1e_2} \ = \ \left(X-\interior(\overline C_1\sqcup \overline C_2) \right) \cup_{\bar g_1\sqcup \bar g_2}(\overline C_1\sqcup \overline C_2)\,. 
$$
Alternatively, $\overline C_1$ and $\overline C_2$ can be joined by an embedded $1$-handle in $X$, the thickening of an arc $\alpha$ in $X-\interior (\overline C_1\sqcup \overline C_2)$ from  $\overline C_1$ to $\overline C_2$.  The result is an embedding $e_1\plus e_2$ of the single cork $(C_1,g_1)\plus(C_2,g_2) = (C_1\plus C_2, g_1\sum g_2)$ (where as noted above the map $g_1\sum g_2$ is only defined up to boundary equivalence unless the $g_i$ are special involutions) 
whose cork twist is independent of $\alpha$.  Indeed, it is readily seen that the single cork twist $\smash{X_{g_1 \sum g_2}^{e_1\plus e_2}}$ is diffeomorphic to the pair of cork twists$X_{g_1g_2}^{e_1e_2}$.  

This process can be iterated to construct the {\it multiple cork twist} $X_{g_1\cdots g_n}^{e_1\cdots e_n}$ of a family $e_1,\dots,e_n$ of disjoint embeddings of corks $(C_1,g_1),\dots,(C_n,g_n)$ in $X$, or a single cork twist $\smash{X_{g_1\sum\cdots\sum g_n}^{e_1\plus \cdots \plus e_n}}$ of an embedding of the boundary sum of the $(C_i,g_i)$'s.  Both twists produce the same smooth $4$-manifold.  This construction will play a key role in what follows.

\vskip -5pt
\vskip -5pt

\subsection*{Trivial Cork Embeddings}

Most explicit corks $(C,g)$ in the literature can be shown to have trivial embeddings in the 4-ball, and thus in every $4$-manifold. In particular, it suffices to prove that the double $C\cup_\id-C$ and twisted double $C\cup_g -C$ are both diffeomorphic to the 4-sphere, often accomplished by an elementary Kirby calculus argument (cf.\ \cite[\S2.6]{akbulut-yasui:stein}).   This is illustrated for the Mazur cork $(\bw,\tau)$ in Figure 2, where the squiggly and straight arrows represent handle slides and cancellations, respectively, and as usual, the $3$ and $4$-handles are not drawn. 

%%%%%%%%%% FIG 2 %%%%%%%%%%
\fig{180}{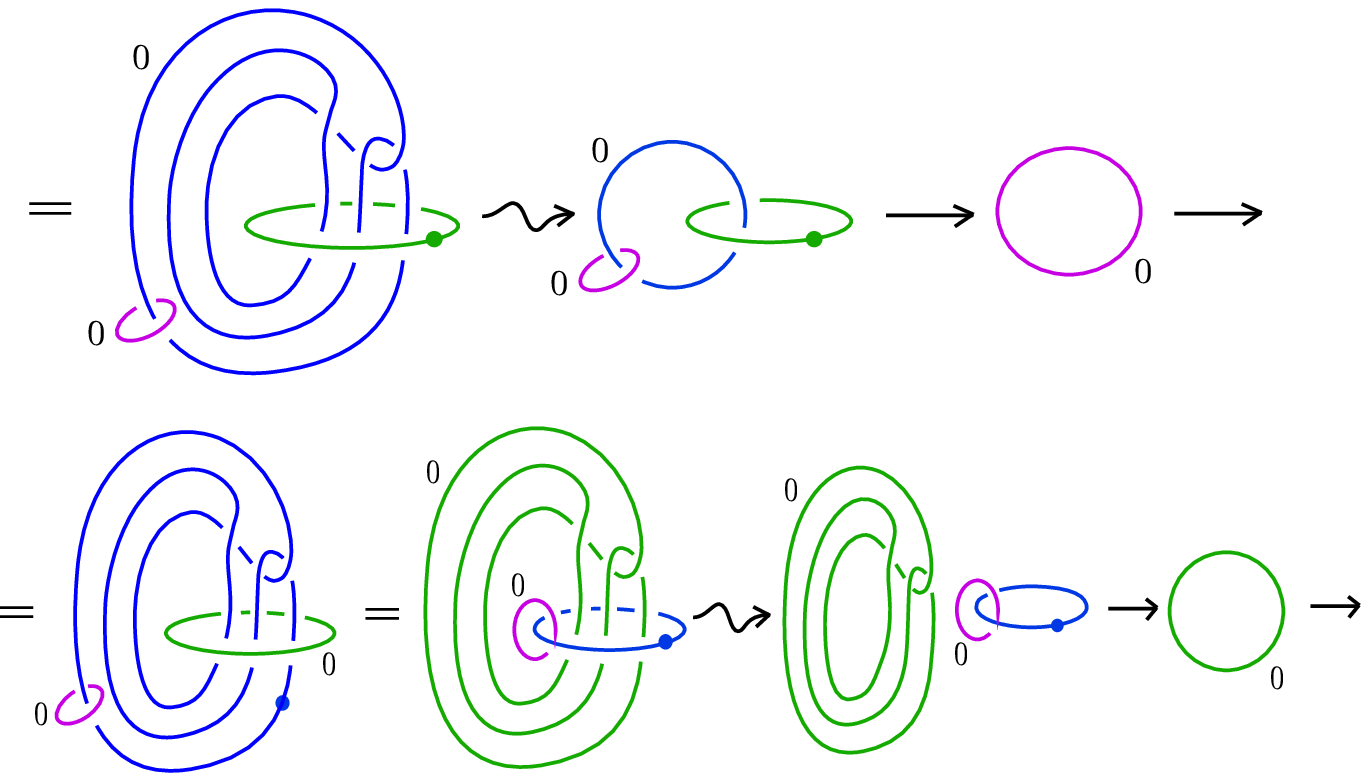}{
  \put(-369,130){$\bw\cup_\id-\bw$}
  \put(-372,36){$\bw\cup_\tau-\bw$}
  \put(-15,130){$B^4\cup B^4$}
  \put(7,36){$B^4\cup B^4$}
\caption{Trivial embedding of the Mazur Cork in $S^4$}
\label{F:MazurCorkTrivial1}}
%%%%%%%%%%%%%%%%%%%%%%%%%%%

\subsection*{Equivariant Corks}   

If $G$ is a subgroup of the diffeomorphism group of $\partial C$ with $(C,g)$ nontrivial for all $g\ne1$ in $G$, then $(C,G)$ is called a {\it $G$-cork}. For cyclic $G$ of finite order $n$, we refer to the corks $(C,g)$ for generators $g$ of $G$ as {\it $n$-corks}.  All explicit corks that have appeared in the literature prior to \cite{tange:cycliccorks} are special $2$-corks; at this time no $\infty$-corks have been shown to exist. 

There is a more general notion, which we call a {\em weakly equivariant cork}, in which the group $G$ is a subgroup of the mapping class group of the boundary, i.e.\ the group of isotopy classes of diffeomorphisms.  In this situation it is more appropriate to use the relation of isotopy, rather than boundary equivalence, because the subgroup of diffeomorphisms of the boundary that extend across the cork need not be normal. Hence the set of boundary equivalent diffeomorphisms does not in general form a group in any natural way.  In the last section, we give a construction of weakly equivariant corks for many groups $G$ that are not subgroups of $\SO(4)$, and in fact that do not act effectively on any homology $3$-sphere.

% Note that there is no obvious way to define a $G$-cork structure on the boundary sum $(C_1,G)\,\plus(C_2,G)$ of two $G$-corks, but there is for special $2$-corks.

In general, if $C$ is a cork with an effective $G$-action on $\partial C$, then an embedding $e:C\hookto X$ will be said to be {\it $G$-effective} if $\smash X_{g_1}^e$ and $X_{g_2}^e$ are smoothly distinct for any $g_1 \ne g_2$ in $G$.  Thus the existence of such an embeddings shows that $(C,G)$ is a $G$-cork.  In this case one has a $G$-action on the set of $4$-manifolds $\{X_{g}^e \ \vert \ g\in G\}$ in the sense that $\smash{(X_{g_1}^e)_{g_2}^{\bar e} = X_{g_1g_2}^e}$
for any two elements $g_1,g_2\in G$, where $\bar e:C\to X_{g_1}^e$ is the obvious embedding induced by $e$. 

For the reader's convenience, we repeat the statement of our main result:

\smallskip

\noindent {\bf Theorem A.} {\it There exist $G$-corks for any finite subgroup $G$ of $SO(4)$.}

\begin{addenda*}
{\bf1)} The proof will show that if $|G|=n$, then the boundary sum $\plus_{n^2}(\bw,\tau)$ of $n^2$ copies of the Mazur cork can be given a $G$-cork structure that has $G$-effective embeddings in any blown-up elliptic surface $E(2k)\cs m\cptwobar$ for $k,m \ge n(n-1)/2$. 

{\bf 2)} More generally, if $\G$ is any finite group that acts effectively on a contractible submanifold of $\R^4$, then essentially the same proof shows that there is a $\G$-cork with an effective embedding into a closed manifold; Theorem \ref{hypthm} can then be used to construct such corks with hyperbolic boundary. 
\end{addenda*}

%%%%%%%%%%%%%%%%%%%%%%%%%%%%%%%%%%%%%%%%%%%%%%%%%%%%%%%%%%%%%%%%%%
%%%%  SECTION 2 : CONSTRUCTION OF EQUIVARIANT CORKS  %%%%%%%%
%%%%%%%%%%%%%%%%%%%%%%%%%%%%%%%%%%%%%%%%%%%%%%%%%%%%%%%%%%%%%%%%%%

\newcommand{\wcheck}[1]{\widecheck{#1}}
\newcommand{\be}{\mathbb E}
\def\E{\be} % this is the blown up Kummer surface
\def\dt{\,\raisebox{-.6ex}{\huge$\hskip-.02in\cdot$\hskip.005in}}
\def\EN{\smash{\E^{\dt}}}     % {\smash{\wcheck\E}}
\def\bs{\mathbb S}
\def\cork{\mathbb T}
\def\corkbar{\overline\cork}

\section{Construction of Equivariant Corks: the proof of \thmref{equiv-cork}}

Our proof of \thmref{equiv-cork} relies on the existence of certain embeddings $e_i$ of the Mazur cork $(\bw,\tau)$ in the blown-up Kummer surface 
$$
\E \ := \ E(2)\cs\cptwobar\,.
$$  
The key input from Seiberg-Witten theory is the count of the number of basic classes in the associated cork twists $\E_\tau^{e_i}$.

\begin{definition}  Let $X$ be a smooth, closed, simply-connected 4-manifold.  If $b_2^+(X)$ is odd and greater than $1$, then $\caln(X)$ will denote the number of Seiberg-Witten basic classes of $X$, and otherwise $\caln(X) = 0$.  For example $\caln(\E) = 2$ (the basic classes are $\pm\cponebar$).
\end{definition}

Akbulut\,\cite{akbulut:contractible} established the nontriviality of $(\bw,\tau)$ by constructing a nontrivial embedding $e_0: \bw\hookto \E$ with {\it reducible} cork twist $\E_\tau^{e_0} \cong 3\cptwo\cs20\cptwobar$, so in particular $\caln(\E_\tau^{e_0})=0$.   It was later observed \cite{bizaca-gompf:elliptic} that such an embedding could be chosen with image in the complement $\EN$ of a {\it nucleus} in $\E$ (see \cite{gompf:nuclei}).  

% Here $E(2)$ is the minimal elliptic surface of Euler characteristic $24$, also known as the Kummer surface (see for example\ \cite{GS:book}).  Nuclei of elliptic surfaces were introduced and first studied by Gompf in \cite{gompf:nuclei}.}  

More recent work of Akbulut and Yasui \cite{akbulut-yasui:knotting-corks} shows that $(\bw,\tau)$ has another nontrivial embedding $e_2:\bw \ \hookto \ \EN$ with $\E_\tau^{e_2}$ {\it irreducible}.  
% Some pictures are drawn in the next section to indicate how $e_2$ is defined.    
The nontriviality of $e_2$ was proved by showing that $\E_\tau^{e_2}$ results from a rational blow-down of $\E$ , 
% \cite{fs:rationalblowdown} 
leaving $\caln$ unchanged, followed by an honest blow-up, doubling $\caln$, so $\caln(\E_\tau^{e_2}) = 4$.
% and so $\caln(\E_\tau^{e_2}) = 4 \ne \caln(\E) = 2$  
(In particular, this follows from Theorem 4.1 for $p=2$, Proposition 5.1 for $n=1$ and $p_1=2$, and Lemma 6.6 in \cite{akbulut-yasui:knotting-corks}.)

As noted in the last section, $(\bw,\tau)$ also embeds {\it trivially} into any $4$-manifold.  Choose one such embedding $e_1:\bw \hookto \EN$.  Thus $e_0,e_1,e_2$ are numbered so that $\caln(\E_\tau^{e_i}) = i\,\caln(\E).$  Only $e_1$ and $e_2$ are needed to prove the following key result, which is a strengthening of an analogous {\it non-compact} embedding theorem of Akbulut and Yasui \cite[Theorem 1.5]{akbulut-yasui:stein}.  

\begin{lemma}\label{L:n-emb}
For each $n>0$, there exists a $2$-cork $(\bs,\sigma)$ that has $n$ disjoint embeddings $s_1,\dots,s_n$ in some closed $4$-manifold $X$, with distinct cork twists $X_\sigma^{s_1} \cong X, X_\sigma^{s_2},\dots,X_\sigma^{s_n}$.  For example the boundary sum $(\bs,\sigma) =\plus_n(\bw,\tau)$ has $n$ such embeddings in the blown-up elliptic surface $X = E(2k)\cs m \cptwobar$ for any $k,m \ge n(n-1)/2$.
\end{lemma}

\proof It suffices to prove the last statement.  First consider the case $k=m=n^2$, and view $X = E(2n^2)\cs n^2\cptwobar$ as the fiber sum of $n^2$ copies of the blown-up Kummer surface $\E = E(2)\cs\cptwobar$ along regular torus fibers in a chosen nucleus.  Denote the copies of $\E$ by $\smash{\E_{ij}}$ for $1\le i,j \le n$.  Choose an embedding $e_{ij}$ of $(\bw,\tau)$ in each summand $\EN_{\!\!\!\!ij}$, with $e_{ij} = e_1$ if $i \le j$ and $e_{ij}=e_2$ if $i > j$.  For $1\le i\le n$, let $s_i$ be the boundary sum $e_{i1}\plus\cdots\plus e_{in}$ of all the embeddings in the ``$i$th row".  Then the $s_i$'s are distinct embeddings of $(\bs,\sigma) = \plus_n(\bw,\tau)$, and can be chosen with disjoint images by choosing the $1$-handles that join the summands to be disjoint.  Furthermore, $s_i$ has $i-1$ nontrivial summands and $n-i+1$ trivial ones, and so $\caln(X_\sigma^{s_i}) = 2^{i-1} \caln(X)$.  Since $\caln(X)\ne0$, the $X_\sigma^{s_i}$ are pairwise distinct.  

Of course one can be more efficient by using only the ``nontrivial" copies of $\E$, i.e.\ $\E_{ij}$ for $i>j$, and putting all the trivial embeddings of the Mazur cork inside one of these.  This handles the smallest case $k=m=n(n-1)/2$, and the fiber sum and blow-up formulas for Seiberg-Witten invariants show that $k$ and $m$ can be increased at will.  \qed 

%\break

\subsection*{Proof of \thmref{equiv-cork}}

Given a finite subgroup $G$ of $\SO(4)$ of order $n$, apply \lemref{n-emb} to produce $n$ disjoint embeddings $s_g$ of a cork $(\bs,\sigma)$ in a closed $4$-manifold $X$, indexed by the elements of $G$, with distinct cork twists $\smash{X_\sigma^{s_g}}$.  Using these cork embeddings, we construct a $G$-cork $(\cork, G)$ and a $G$-effective embedding $t:\cork\to X$, as follows.

The underlying contractible manifold $\cork$ is the boundary sum $\natural_{\,n}\,\bs$ of $n$ copies of $\bs$.  To define the $G$ action on $\partial\cork$, it is convenient to present $\cork$ as a cork twist on a diffeomorphic copy $\corkbar$ of $\cork$ that supports a natural $G$-action, namely the equivariant boundary sum
$$
\corkbar \ = \ B^4 \ \plus \ (G\times\bs)\,.
$$
Here $G$ acts on $B^4$ linearly, and on $G\times\bs$ by left multiplication on the first factor and trivially on the second.  

Before defining $\cork$, we show how to use the embeddings $s_g$ of $\bs$ to construct an embedding 
$$
\bar t:\corkbar\hookto X.
$$
Start with the disjoint union of the $4$-ball with its linear $G$-action, and $n$ copies $\bs_g$ of $\bs$, indexed by the elements of $G$.  To get $\corkbar$, add $1$-handles joining $b_g\in\partial B$ to $x_g\in \partial\bs_g$, where $\{b_g\st g\in G\}$ is a principal $G$-orbit in $\partial B$ and the $x_g \in \partial\bs_g$ correspond to a chosen $x\in\partial\bs$.  Then $G$ acts linearly on $B$ and permutes the copies of $\bs_g$ by left multiplication on the subscript.  Now identify $\bs_g$ with the image $s_g(\bs)$ in $X$, and $B^4$ with a small $4$-ball $B$ disjoint from the $\bs_g$'s.  Joining $B$ by embedded $1$-handles to the $\bs_g$'s gives the desired embedding $\bar t$.  

To obtain $\cork$, we twist a shrunken copy of the cork $1\times\bs$ in $\corkbar$.  To make this precise, recall that $\corkbar$ contains $n$ copies $\bs_g = g\times\bs$ of $\bs$, the images of the embeddings $e_g:\bs\hookto\corkbar$ sending $x$ to $(g,x)$.  Consider an embedding $s:\bs\hookto\bs$ that shrinks $\bs$ inside itself, that is, $s$ is the identity off of a boundary collar $\partial\bs\times [0,1)$, and maps $(x,t)$ to $(x,(t+1)/2)$ inside the collar.  Then $e = e_1\circ s$ embeds $\bs$ onto a shrunken copy of $\bs_1$.  We define $\cork$ to be the cork twist associated with this embedding,
$$
\cork \ = \ \corkbar_\sigma^{\,e}.
$$ 
Since the $\partial\cork = \partial\corkbar$, there is still a $G$-action on $\partial\cork$, and this defines our cork $(\cork,G)$.  Note that $\cork$ is actually diffeomorphic to $\corkbar$, and thus to $\natural_{\,n}\,\bs$, since $\natural$ is a well defined operation, but for our purposes it is most convenient to describe $\cork$ as a cork twist of $\corkbar$.    

Now observe that the embedding $\bar t:\corkbar\hookto X$ above induces an embedding
$$
t:\cork \ \hookto \ X_\sigma^{s_1}
$$
since $\cork =\corkbar_\sigma^e$.  Furthermore, twisting this embedding of $\cork$ by an element $g\in G$ just transfers the cork twist from $\bs_1$ to $\bs_g$, that is 
$$  
(X_\sigma^{s_1})_g^t \ = \ X_\sigma^{s_g}.
$$
Since the smooth $4$-manifolds $X_\sigma^{s_g}$ are distinct for $g\in G$, this shows that $t$ is a $G$-effective embedding, and so $(\cork,G)$ is a $G$-cork.  This completes the proof of Theorem 1.1.   \qed

\begin{remark*}
Even in the case $\G=\Z_2$ this result can give something new. Applying the construction from \thmref{equiv-cork} to the free $\Z_2$ action on $S^3$ extended across $B^4$ we get a $2$-cork with free action on the boundary.
\end{remark*}

\subsection*{Proof of the Addenda to \thmref{equiv-cork}}
The first addendum to the theorem follows from this proof by using $(\bs,\sigma) = \plus_n(\bw,\tau)$ and $X = E(2k)\cs m \cptwobar$, as provided by the lemma.  Note that in the proof, $X_\sigma^{s_1}$ is diffeomorphic to $X$ since $s_1$ is a trivial cork embedding, so $t$ can be viewed as an embedding of $\natural_{\,n^2}\,\bw \hookto X$.

With regard to the second addendum, if a finite group $\G$ acts on a compact contractible submanifold of $\R^4$, we may repeat the argument replacing $B^4$ by the contractible submanifold to produce a $G$-cork $\cork$. To build a $\G$-cork with hyperbolic boundary, let $\bu$ be an invertible cobordism from $\partial\cork$ to a hyperbolic $3$-manifold $M$ with inverse $\bv$ as given by Theorem \ref{hypthm}. Then  
$$
\cork\cup_{\partial\cork} \bu \ \subset \ \cork\cup_{\partial \cork} \bu \cup_M \bv\cong \cork
$$
and $\cork\cup_{\partial \cork}\bu$ inherits a $\G$ action so twisting it via $g$ has the same effect as  twisting $\cork$ since $g$ extends across $\bv$.      
      
\begin{remark*}
From the construction, we see that our $G$-corks are boundary-connected sums of Stein manifolds, and hence are Stein.  In contrast to the argument in~\cite{tange:cycliccorks}, this fact does not play any role in our verification that our corks are effective.
\end{remark*}

%%%%%%%%%%%%%%%%%%%%%%%%%%%%%%%%%%%%%%%%%%%%%%%%%%%%%%%%%%%%%%%%%%
%%%%  SECTION 3 : WEAKLY EQUIVARIANT CORKS  %%%%%%%%
%%%%%%%%%%%%%%%%%%%%%%%%%%%%%%%%%%%%%%%%%%%%%%%%%%%%%%%%%%%%%%%%%%

\section{Weakly equivariant corks}\label{S:w-eq}
%The question of determining whether finite subgroups of the mapping class group of a manifold of dimension lift to the diffeomorphism group is known as the Nielsen realization problem. Such lifts were shown to exist in dimension $2$, but in higher dimension, they typically do not 
%exist~\cite{raymond-scott:nielsen}. 
In this section, we construct examples of weakly equivariant corks for certain finite groups that are not subgroups of $\SO(4)$. In fact, these groups cannot act on any homology sphere, so there are no corresponding equivariant corks. This will prove \thmref{w-eq} in the introduction.

Fix $n$, and let $G$ be the elementary abelian group $(\Z_2)^n$.  It is known that for $n>3$, this group does not to act effectively on any homology $3$-sphere~\cite[Proposition 3]{zimmermann:classification}.  In the following paragraphs, we show how to construct a weak $G$-cork $\bv$.  

Apply \lemref{n-emb} to get a $2$-cork $(\bs,\sigma)$ with $2^n$ inequivalent embeddings $s_1,\dots,s_{2^n}$ in some $4$-manifold $X$, meaning their cork twists $X_\sigma^{s_i}$ are distinct smooth 4-manifolds, with $X_\sigma^{s_1}\cong X$.  For example, $\bs$ could be the boundary sum of $2^n$ Mazur corks, with $X = E(2^{2n+1})\cs 2^{2n}\cptwobar$.

\def\bvbar{\overline\bv}

As in the proof of \thmref{equiv-cork}, the cork $\bv$ will be defined as a suitable cork twist of $\bvbar$, a diffeomorphic copy of $\bv$.  To define $\bvbar$, consider a full binary tree $T$ of height $n$, built from the bottom up, as shown in \figref{BinaryCork} for the case $n=4$.  Thus $T$ has one vertex at the root, two at the first level, four at the second level, etc.  At the top there are $2^n$ vertices, the {\it leaves} of $T$.  To get $\bvbar$, replace the black dots by $4$-balls, the white dots by copies of the cork $\bs$ (referred to as the {\it leaves} of the cork), and the edges by $1$-handles.  Also choose an {\it equatorial $2$-sphere} on the boundary of each $4$-ball that separates the attaching foot of the 1-handle from below (if any) from the attaching feet of the 1-handles above. 

%%%%%%%%%% FIG 1 %%%%%%%%%%
\fig{100}{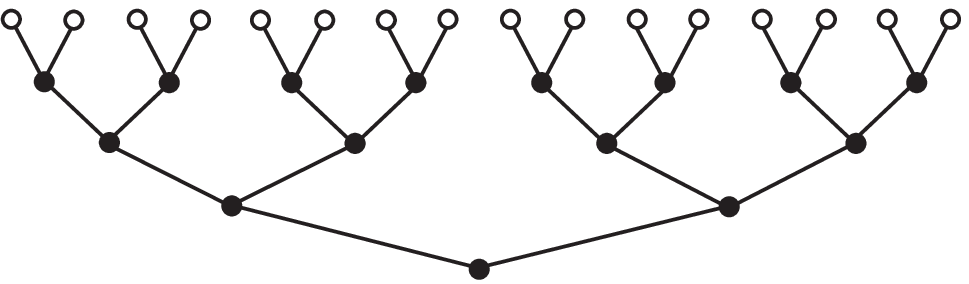}{
\caption{A weak $(\bz_2)^4$-cork}
\label{F:BinaryCork}}
%%%%%%%%%%%%%%%%%%%%%%%%%%%

Let $\tau_0,\dots,\tau_{n-1}$ denote the generators of the $\Z_2$-factors in $G = (\Z_2)^n$, and let $\tau_k$ act on $\partial\bvbar$ by performing half-Dehn twists on all of the $2^k$ equatorial $2$-spheres at level $k$.  
%Note that $\tau_0$ reverses the order of the leaves at the top, $\tau_1$ independently reverses the orders of the first and second halves of the leaves, etc.
Then it is readily verified that $\tau_k$ is of order $2$ in the mapping class group $\overline\calm$ of $\partial\bvbar$ (this is clear for $\tau_0$, and in general $\tau_k^2$ is isotopic to $\tau_{k-1}^2$) and that the action of the $\tau_k$'s extends to an embedding of $G$ in $\overline\calm$.

Now define $\bv$ to be the cork twist of $\bvbar$ along (a shrunken copy of) the leftmost leaf in $\bvbar$.  Then $\partial\bv = \partial\bvbar$, and so there is an induced embedding of $G$ in the mapping class group of $\partial\bv$.  To see that this defines a weak $G$-cork structure on $\bv$, just choose an embedding $e:\bv\hookto X$ that restricts to the embeddings $s_1,\dots,s_{2^n}$ on the leaves of $\bv$ (from left to right) and then it is clear that the cork twists $X_g^e$ and $X_h^e$ are not diffeomorphic for distinct elements $g$ and $h$ in $G$. This completes the proof of \thmref{w-eq}.  \qed

%\newpage

\def\cprime{$'$}
\providecommand{\bysame}{\leavevmode\hbox to3em{\hrulefill}\thinspace}
\providecommand{\KRhref}[2]{%
  \href{http://www.ams.org/mathscinet-getitem?mr=#1}{#2}
}
\providecommand{\href}[2]{#2}

\end{document}